\documentclass[12pt, reqno]{amsart}
\usepackage{amsmath, amsthm, amscd, amsfonts, amssymb, graphicx, color}
\usepackage[bookmarksnumbered, colorlinks, plainpages]{hyperref}

\usepackage[all]{xy}
\newtheorem{theorem}{Theorem}[section]

\newtheorem{proposition}[theorem]{Proposition}
\newtheorem{corollary}[theorem]{Corollary}
\theoremstyle{definition}
\newtheorem{definition}[theorem]{Definition}
\newtheorem{example}[theorem]{Example}

\theoremstyle{remark}
\newtheorem{remark}[theorem]{Remark}
\numberwithin{equation}{section}

\begin{document}
	\setcounter{page}{1}
	\title[\textit{DW}-compact operators on Banach lattices]{\textit{DW}-compact operators on Banach lattices}
	
    \author[J.X. Chen]
       {Jin Xi Chen}
       \address{School of Mathematics, Southwest Jiaotong University, Chengdu 610031, China}
       \email{jinxichen@swjtu.edu.cn}

	 \author[J. Feng]
	  {Jingge Feng}
	 \address { School of Mathematics, Southwest Jiaotong University, Chengdu 610031, China}
	
\email{fjg0825@my.swjtu.edu.cn}

	

	\begin{abstract}
   This paper is devoted to the study of \textit{DW}-compact operators, that is, those  operators  which map disjointly weakly compact sets in a Banach lattice onto relatively compact sets. We show that  \textit{DW}-compact operators are precisely the operators which are both Dunford-Pettis and  \textit{AM}-compact. As an application, Banach lattices with the property that every disjointly weakly compact set is a limited (resp. Dunford-Pettis) set, are characterized by using \textit{DW}-compact operators.

	\end{abstract}
\maketitle\baselineskip 4.80mm
	
\section{Introduction and preliminaries}
Throughout this paper, $X, Y$ will denote real Banach spaces, and
	$E, F$ will denote real Banach lattices. $B_X$ is the closed unit ball of $X$, and $E^+$ is the positive cone of the Banach lattice $E$.  The solid hull of a subset $A$ of $E$ is denoted by $Sol(A):=\lbrace y\in E:|y| \leq |x|$ \ for \ some\ $x\in A\rbrace$.

\par Following W. Wnuk \cite{WOrder}, a bounded subset $A$ of the Banach lattice $E$ is called a \textit{disjointly weakly compact set}
if every disjoint sequence from $Sol(A)$ converges weakly to zero. Weakly precompact sets and order bounded sets in Banach lattices are all disjointly weakly compact \cite[Proposition 2.5.12 iii)]{Meyer}. It is well known that the dual $E^{\,\prime}$ has order continuous norm if and only if $B_E$ is disjointly weakly compact. Recently, the disjointly weak compactness in Banach lattices has been extensively investigated by Xiang, Chen and Li in \cite{disjointly weakly compact set}. Therein, they proved that every disjointly weakly compact subset of $E$ is relatively weakly compact if and only if $E$ is a $KB$-space \cite[Proposition 2.6]{disjointly weakly compact set}. Naturally, a class of operators called disjointly weakly compact operators was introduced. An operator $T\in\mathcal{L}(X,E)$ from a Banach space into a Banach lattice $E$ is called a \textit{disjointly weakly compact operator} if $TB_X$ is a disjointly weakly compact subset of $E$.  $T\in\mathcal{L}(X,E)$ is disjointly weakly compact if and only if $T^{\prime}$ is an order weakly compact operator \cite[Theorem 3.3]{disjointly weakly compact set}.

\par Besides compact operators and Dunford-Pettis operators, \textit{AM}-compact operators introduced by P. G. Dodds and D. H. Fremlin \cite{Dodds and Fremlin PL-compact},  are also a class of operators which combine the order structure and the compactness property. Let us recall that an operator $T\in\mathcal{L}(E,X)$ from a Banach lattice into a Banach space  is called an \textit{AM-compact operator} if $T$ carries each order bounded subset of $E$ onto a relatively compact set in $X$. The name \textit{AM-compact operator} is justified by a simple observation:  If $E$ is an \textit{AM}-space with order unit, then \textit{AM}-compact operators are precisely those operators which are compact. From a result of A. Grothendieck it follows that $T\in\mathcal{L}(E,X)$ is \textit{AM}-compact if and only if $T^{\prime}B_{X^\prime}$ is $|\sigma|(E^\prime,E)$--totally bounded. See, e.g., \cite[Theorem 3.27, Theorem 3.55; Exercise 8, p.180 ]{Positive}. Theorem 4.7 of \cite{Dodds and Fremlin PL-compact} asserts that  the set of all regular \textit{AM}-compact operators from $E$ to $F$ forms a band in $\mathcal{L}^{r}(E,F)$ when the range space $F$ is a Banach lattice with order continuous norm.

 \par Recall that a bounded subset $A$ of  $E$ is called an \textit{$L$-weakly compact set} if  $\|x_{n}\|\to 0$ for each disjoint sequence $(x_n)$ in $Sol(A)$. Every $L$-weakly compact set is  relatively weakly compact, and $L$-weakly compact sets and  relatively weakly compact sets coincide in a Banach lattice with the positive Schur property \cite[pp.212--213]{Meyer}. Here, we say a Banach lattice has the \textit{positive Schur property} if every disjoint weakly null sequence in $E$ is norm null, or equivalently, if every weakly null sequence of $E^+$ is norm null. Abstract $L$-spaces are typical examples of Banach lattices with the positive Schur property. In fact, the positive Schur property of Banach lattices can also be characterized as the property that   every disjointly weakly compact set  is  $L$-weakly compact \cite[Corollary 2.7]{disjointly weakly compact set}. Recently, Hajji and  Mahfoudhi \cite{LW1} introduced a class of operators that they called \textit{$LW$-compact operators}, that is, those operators which map $L$-weakly compact sets onto relatively compact sets.

\par Observe that  weakly precompact sets, $L$-weakly compact sets  and order bounded sets in Banach lattices are all disjointly weakly compact sets.  Motivated by this observation, we introduce a class of operators which map disjointly weakly compact sets onto relatively compact sets. We call a bounded linear operator $T:E\to X$ \textit{\textit{DW}-compact} if $T(A)$ is  relatively compact for each disjointly weakly compact subset $A$ of $E$. Clearly, a \textit{DW}-compact operator is necessarily Dunford-Pettis, \textit{AM}-compact, and \textit{LW}-compact.  Indeed, we will see that \textit{DW}-compact operators are precisely those operators which are both Dunford-Pettis and  \textit{AM}-compact (Theorem \ref{AM-compact and almost DP}). Some results on the domination problem for positive \textit{DW}-compact operators are  obtained. As applications, we also give some characterizations of order-topological properties (e.g., the positive Schur property, the (d)-DP property and the (d)-DP$^*$ property) of Banach lattices in terms of \textit{DW}-compact operators.

\par The definitions and notions from Banach lattice theory which appear here, are standard. We refer the reader to the references \cite{Positive, Meyer}. For locally convex-solid topologies on vector lattices see \cite{Locally Solid, Positive}.

\section{\textit{DW}-compact operators on Banach lattices}

We start with several results on the properties of disjointly weakly compact sets. By definition, a \textit{DW}-compact operator  maps disjointly weakly compact sets onto relatively compact ones. The following result shows that disjointly weakly compact sets can also be identified by \textit{DW}-compact operators with ranges in $c_0$.

\begin{theorem}\label{disjointly weakly compact by DW-compact}
For a bounded subset $A$ of $E$ the following statements are equivalent.
\begin{enumerate}
  \item $A$ is disjointly weakly compact.
  \item Each \textit{DW}-compact operator from $E$ into an arbitrary Banach space carries $A$ onto a relatively compact set.
  \item Each \textit{DW}-compact operator from $E$ into $c_0$ carries $A$ onto a relatively compact set.
\end{enumerate}
\end{theorem}
\begin{proof}
$(1)\Rightarrow (2)$ and $(2)\Rightarrow (3)$ are obvious.
\par $(3)\Rightarrow (1)$ In view of Theorem 2.3 of \cite{disjointly weakly compact set}, it suffices to prove that every order bounded disjoint sequence of $E^\prime$ converges uniformly to zero on $A$. Now, let $(f_n)_{n=1}^{\infty}$ be an order bounded disjoint sequence in $E^\prime$. Note that $f_{n}\xrightarrow{w}0$ and consider the operator $T:E\to c_0$ defined by $$Tx=(f_{1}(x),f_{2}(x),\cdots,f_{n}(x),\cdots),\,\,\,x\in E.$$We claim that $T$ is \textit{DW}-compact. To this end, let $D$ be a disjointly weakly compact subset of $E$. From Theorem 2.3 of \cite{disjointly weakly compact set} it follows that $(f_n)_{n=1}^{\infty}$ converges to zero uniformly on $D$, that is, $\sup_{x\in D}|f_{n}(x)|\to{0}$. This implies that $T(D)\subset c_0$ is relatively compact. Hence, $T$ is \textit{DW}-compact. Now, by our hypothesis, $T(A)$ is a relatively compact subset of $c_0$, equivalently, $\sup_{x\in A}|f_{n}(x)|\to{0}$. That is,  $(f_n)_{n=1}^{\infty}$ converges uniformly to zero on $A$. Again, by  \cite[Theorem 2.3]{disjointly weakly compact set} we can see that $A$ is disjointly weakly compact.
\end{proof}

Since  $B_E$ is disjointly weakly compact if and only if $E^\prime$ has order continuous norm, the following corollary is an immediate consequence of Theorem \ref{disjointly weakly compact by DW-compact} and tells us when each \textit{DW}-compact operator is compact.

\begin{corollary}\label{DW-compact is compact}
For a Banach lattice $E$ the following statements are equivalent.
\begin{enumerate}
  \item Every \textit{DW}-compact operator from $E$ into an arbitrary Banach space is a compact operator.
  \item Every \textit{DW}-compact operator from $E$ into $c_0$  is a compact operator.
  \item $E^{\,\prime}$ has order continuous norm.
  \end{enumerate}
\end{corollary}

The disjointly weak compactness of a set can be preserved by order bounded operators.
\begin{proposition}\label{preservation of disjointly weakly compact}
Let $T:E\to F$ be an order bounded operator between Banach lattices and let $A\subset E$ be a disjointly weakly compact set. Then $T(A)$ is likewise disjointly weakly compact.
\end{proposition}

\begin{proof}
We assume without loss of generality that $A$ is solid. To prove that $T(A)$ is disjointly weakly compact, by \cite[Theorem 2.3]{disjointly weakly compact set} we have to show that every order bounded sequence $(f_n)\subset (F^\prime)^+$ satisfying $f_{n}\xrightarrow{w^*}0$ converges  uniformly to zero on $T(A)$. First, we claim that $|T^{\prime}f_{n}|\xrightarrow{w^*}0$ in $E^\prime$. Indeed, this follows from the equalities
  $$|T^{\prime}f_{n}|(x) =\sup_{u\in [-x,x]} |f_{n}(Tu)|\xrightarrow{n\to\infty} 0,\,\,\,x\in E^+ $$since $T[-x,x]$ is order bounded and hence is disjointly weakly compact. Therefore, $(|T^{\prime}f_{n}|)_{n=1}^{\infty}$ is an order bounded and weak$^*$-null sequence in $E^\prime$ since $T^\prime$ is likewise order bounded. The disjointly weak compactness of $A$ implies that
    $$\sup_{x\in A}|f_{n}(Tx)|\leq\sup_{x\in A}|T^{\prime}f_{n}|(|x|)=\sup_{x\in A}\left||T^{\prime}f_{n}|(x)\right|\to{0}.$$That is, $T(A)$ is a disjointly weakly compact subset of $F$.
\end{proof}

 We can characterize \textit{DW}-compact operators in terms of compact and disjointly weakly compact operators.

 \begin{theorem}\label{disjointly and DW}
	For a bounded linear operator $T:E\rightarrow X$ from a Banach lattice to a Banach space the following statements are equivalent:
\begin{enumerate}
	\item $T$ is a \textit{DW}-compact operator.
	
	\item For every Banach space $Y$ and every disjointly weakly compact operator $S:Y\rightarrow E$, the composite $TS$ is a compact operator.
	
	\item For every disjointly weakly compact operator $S:\ell^{1}\rightarrow E$, the composite $TS$ is a compact operator.
\end{enumerate}
\end{theorem}

\begin{proof}
 It suffices  to prove (3)$\Rightarrow$(1). Let $A$ be a disjointly weakly compact subset of $E$  and let $(x_n)$ be a sequence in $A$. We have to show that $(Tx_n)$ has a norm convergent subsequence. Consider the operator $S:\ell^1\rightarrow E$  defined by
 $$S((\lambda_n))=\sum\limits_{n=1}^{\infty} \lambda_{n}x_{n}, \,\,\, (\lambda_n)\in \ell^1 $$Then, $S$ is a disjointly weakly compact operator since $SB_{\ell_1}\subset \overline{\{x_{n}:n\in \mathbb{N}\}^{bc}}$, where $\overline{\{x_{n}:n\in \mathbb{N}\}^{bc}}$, the closed absolutely convex hull of $(x_n)$, is disjointly weakly compact (see \cite[Theorem 2.3]{disjointly weakly compact set}). Thus, by our hypothesis the composite  $TS:\ell^1\rightarrow X$ is a compact operator. It follows that $(Tx_n)=(TSe_n)$ has a  convergent subsequence. This implies that $T$ is \textit{DW}-compact.
\end{proof}

Let us recall that a  bounded linear operator $T: E\rightarrow X$ from a Banach lattice to a Banach space  is called an \textit{almost Dunford-Pettis operator} if $\|Tx_{n}\| \rightarrow 0$ for each disjoint weakly null sequence $(x_n)$ in $E$, or equivalently, if $\|Tx_{n}\| \rightarrow 0$ whenever $0\leq x_{n}\xrightarrow{w}0$ in $E$ (\cite{some,Sanchez}). The identity operator $I:L^{1}[a,b]\to L^{1}[a,b]$ is an example of an almost Dunford-Pettis which is not \textit{DW}-compact. On the other hand, the identity $I:c_0\to c_0$ is an \textit{AM}-compact operator which is not \textit{DW}-compact.  The following result shows that a \textit{DW}-compact operator is indeed the combination of an \textit{AM}-compact operator and an (almost) Dunford-Pettis operator.
\begin{theorem}\label{AM-compact and almost DP}
A bounded linear operator $T:E\to X$ is  \textit{DW}-compact  if and only if $T$ is both \textit{AM}-compact and (almost) Dunford-Pettis.
\end{theorem}

\begin{proof}
If $T$ is a \textit{DW}-compact operator, then it is obvious that $T$ is both \textit{AM}-compact and  Dunford-Pettis.
\par For the converse, let $T:E\to X$ be both \textit{AM}-compact and almost Dunford-Pettis, and let $A$ be a disjointly weakly compact subset of $E$. By definition, we can assume without loss of generality that $A$ is solid. For every disjoint sequence $(x_n)\subset A$, we have $\|Tx_n\|\to 0$ since $T$ is almost Dunford-Pettis. Then in view of Theorem 4.36 of \cite{Positive}, for each $\varepsilon>0$  there exists some $u\in{E^{+}}$ lying in the ideal generated by $A$ such that
	$\|T[(|x|-u)^{+}]\|<\frac{\varepsilon}{2}$ holds for all $x\in A$. Therefore, from the identity $|x|=|x|\wedge u+(|x|-u)^+$ we can see that
$$T(A)\subset T[-u,u]+\varepsilon B_{X}$$
Since $T$ is \textit{AM}-compact, $T[-u,u]$ is relatively compact. It follows that  $T(A)$ is relatively compact. This implies that $T$ is \textit{DW}-compact.
\end{proof}
As a byproduct of Theorem \ref{AM-compact and almost DP},  an \textit{AM}-compact operator $T:E\to X$ is Duford-Pettis if and only if it is almost Dunford-Pettis. By Theorem \ref{AM-compact and almost DP} we can easily see that the identity operator $I:c_0\to c_0$ is not \textit{DW}-compact while its adjoint is a \textit{DW}-compact operator. On the other hand, the identity operator $I:\ell^1\to \ell^1$ is an example of a \textit{DW}-compact operator whose adjoint is not \textit{DW}-compact.

\par As we have already seen, the identity operator $I:L^{1}[a,b]\to L^{1}[a,b]$ is not \textit{DW}-compact even though $L^{1}[a,b]$ has the positive Schur property. However, this situation does not happen to a Banach lattice with the Schur property.

\begin{theorem}\label{Each linear operator is DW}  For a Banach lattice $E$ the following assertions are equivalent:
\begin{enumerate}

  \item $E$ has the Schur property.
  \item Each disjointly weakly compact subset of $E$ is relatively compact, that is, the identity operator $I:E\to E$ is \textit{DW}-compact.
  \item Each bounded linear operator $T:E\to X$ from $E$ to an arbitrary Banach space $X$ is \textit{DW}-compact.
  \item Each bounded linear operator $T:E\to\ell_\infty$  is \textit{DW}-compact.

\end{enumerate}
\end{theorem}
\begin{proof}
$(1)\Rightarrow (2)$ Let $A$ be a  disjointly weakly compact subset of $E$.  Since $E$ has the Schur property, $E$ is a \textit{KB}-space. From Proposition 2.6 of \cite{disjointly weakly compact set} it follows that $A$ is relatively weakly compact. Once again, the Schur property of $E$ implies that $A$ is relatively compact.
\par $(2)\Rightarrow (3)$ and $(3)\Rightarrow (4)$ are obvious.
\par $(4)\Rightarrow (1)$ We assume by way of contradiction that $E$ does not have the Schur property. Then there exists a sequence $(x_n)\subset E$ such that $x_{n}\xrightarrow{w}0$ and $\|x_n\|=1$ for each $n\in\mathbb{N}$. Choose $f_{n}\in E^\prime$ satisfying $\|f_n\|=1$ and $f_n(x_n)=1$ \,\, $(n\in\mathbb{N})$. Let us define a bounded linear operator $T:E\to\ell_\infty$ by $Tx=(f_{n}(x))_{n=1}^{\infty}$. Then $\|Tx_n\|\geq|f_{n}(x_n)|=1$. This implies that $T$ is not Dunford-Pettis and hence not \textit{DW}-compact.
\end{proof}
It should be noted that a Dunford-Pettis operator is not necessarily \textit{DW}-compact.
\begin{example}\label{DP is not DW}
   We know that every  operator $T\in \mathcal{L}(\ell^{\infty},c_0)$ is weakly compact and Dunford-Pettis since $\ell^{\infty}$ is a Grothendieck space with the Dunford-Pettis property. However, $\mathcal{L}(\ell^{\infty},c_0)\neq \mathcal{K}(\ell^{\infty},c_0)$. This implies that there exists an operator $S\in \mathcal{L}(\ell^{\infty},c_0)$ which is not \textit{DW}-compact since $B_{\ell^{\infty}}$ is disjointly weakly compact.
\end{example}

 A bounded subset $A$ of  $X$ is said to be \textit{weakly  precompact} or \textit{conditionally weakly compact} if every sequence in $A$ contains a weak Cauchy subsequence. Indeed,  Based on the work of Odell and Stegall, Rosenthal \cite[p.377]{1977Point} pointed out that a subset $A$ of  $X$ is weakly precompact if and only if for every Banach space $Y$ and every Dunford-Pettis operator $T:X\to{Y}$, $T(A)$ is  relatively compact.  Ghenciu \cite{Ioana} showed that $A\subset X$ is weakly precompact if and only if $T(A)$ is relatively compact for every Dunford-Pettis operator $T:X\to{c_0}$ \cite[Theorem 1 \& Corollary 9]{Ioana}. This  implies that Dunford-Pettis operators with ranges in $c_0$ can also be employed to identify weakly precompact sets.  Xiang, Chen and Li \cite{XCL} recently considered weak precompactness in Banach lattices. They proved that the solid hull of every weakly precompact subset of $E$ is likewise weakly precompact if and only if every disjointly weakly compact subset of $E$ is weakly precompact, or equivalently, if  every order interval in $E$ is weakly precompact \cite[Theorem 2.3 \& 2.4]{XCL}. In particular, if $E$ has order continuous norm or  $E^\prime$ has a weak order unit, then every weakly precompact subset of $E$ has a weakly precompact solid hull \cite[Proposition 2.5.12 ii)] {Meyer}. Combine the above comments and Theorem \ref{disjointly weakly compact by DW-compact} we can see the following result is evident.

\begin{proposition}\label{DP is DW}
For a Banach lattice $E$ the following assertions are equivalent:
\begin{enumerate}
  \item Every disjointly weakly compact subset of $E$ is weakly precompact, i.e., every order interval in $E$ is weakly precompact.
  \item Every Dunford-Pettis operator $T:E\to X$ from $E$ to an arbitrary Banach space $X$ is \textit{DW}-compact.
  \item Every Dunford-Pettis operator $T:E\to c_0$  is \textit{DW}-compact.
\end{enumerate}
\end{proposition}
One may ask when an almost Dunford-Pettis operator is \textit{DW}-compact. To answer this we have to mention totally bounded sets with respect to the absolute weak topologies on Banach lattices. Recall that the absolute weak topology $|\sigma|(E,E^\prime)$ on a Banach lattice $E$ is the locally convex-solid topology generated by the family of lattice seminorms $\{\rho_{x^\prime}:x^\prime \in E^\prime\}$, where $\rho_{x^\prime}$ is defined by $\rho_{x^\prime}(x)=|x^{\prime}|(|x|)$ for $x\in E$.  A subset $A$ of $E$ is $|\sigma|(E,E^\prime)$--totally bounded if for every $\varepsilon>0$ and every finite collection $\{x_{1}^{\prime}, x_{2}^{\prime},\cdots,x_{n}^{\prime}\}\subset E^\prime$ there exists a finite subset $\Phi$ of $A$ such that $A\subset \Phi+\bigcap_{i=1}^{n}\{x:|x_{i}^{\prime}|(|x|)<\varepsilon\}$. Recently,  Ardankani and  Chen \cite{totally bounded in absolute} showed that a  subset $A$ of  $E$ is $|\sigma|(E,E^\prime)$-totally bounded if and only if  $T(A)\subset c_0$ is relatively compact for every almost Dunford-Pettis operator $T:E\to{c_0}$. They also proved that every disjointly weakly compact subset of $E$ is $|\sigma|(E,E^{\prime})$-totally bounded if and only if every order interval in $E$ is $|\sigma|(E,E^{\prime})$-totally bounded, or equivalently, $E^\prime$ is discrete \cite[Corollary 2.13\&Remark 2.14]{totally bounded in absolute}.
\begin{proposition}\label{almost DP is DW}
For a Banach lattice $E$ the following assertions are equivalent:
\begin{enumerate}
  \item Every disjointly weakly compact subset of $E$ is $|\sigma|(E,E^{\prime})$-totally bounded, i.e.,  $E^\prime$ is discrete.
  \item Every almost Dunford-Pettis operator $T:E\to X$ from $E$ to an arbitrary Banach space $X$ is \textit{DW}-compact.
  \item Every almost Dunford-Pettis operator $T:E\to c_0$  is \textit{DW}-compact.
\end{enumerate}
\end{proposition}

\begin{proof}
 $(1)\Rightarrow(2)$  Assume that every disjointly weakly compact subset of $E$ is $|\sigma|(E,E^{\prime})$-totally bounded. Then this implies that every weakly precompact subset of $E$ is $|\sigma|(E,E^{\prime})$-totally bounded. From Theorem 2.6 of \cite{totally bounded in absolute} it follows that every almost Dunford-Pettis operator $T:E\to X$ from $E$ to an arbitrary Banach space $X$
is Dunford-Pettis. On the other hand, since every disjointly weakly compact subset of $E$ is $|\sigma|(E,E^{\prime})$-totally bounded and hence is weakly precompact, by  Proposition \ref{DP is DW}  every almost Dunford-Pettis operator $T:E\to X$  is \textit{DW}-compact.
\par $(2)\Rightarrow(3)$ Obvious.
\par $(3)\Rightarrow (1)$ Let $A$ be a disjointly weakly compact subset of $E$. For each almost Dunford-Pettis operator $T:E\to c_0$, by our hypothesis $T$ is \textit{DW}-compact. Hence, $TA$ is relatively compact in $c_0$. This implies that $A$ is $|\sigma|(E,E^{\prime})$-totally bounded \cite[Theorem 2.10]{totally bounded in absolute}.
\end{proof}

Recall that an operator $T\in \mathcal{L}(E,X)$  is \textit{\textit{LW}-compact}  if $T$ maps \textit{L}--weakly compact subsets of $E$ onto  relatively  compact subsets of $X$ \cite{LW1}. Every \textit{DW}-compact operator is obviously \textit{LW}-compact.  We can easily see that the converse does not necessarily hold. For instance, if $E=c_{0}$ or $\ell^{p} (1<p\leq\infty)$, then every $L$--weakly compact subset of $E$ is relatively compact. Hence, the identity operator $I:E\to E$ is $LW$-compact. However, $I$ is not \textit{DW}-compact since $B_E$ is disjointly weakly compact.

\begin{theorem}\label{LW is DW}
	For  a Banach lattice $E$  the following statements  are equivalent.
\begin{enumerate}	
	\item $E$ has the positive Schur property.
	
	\item Every \textit{LW}-compact operator $T:E\rightarrow X$ from $E$ into an arbitrary Banach space $X$ is \textit{DW}-compact.
	
	\item Every \textit{LW}-compact operator $T:E\rightarrow \ell^\infty$ is \textit{DW}-compact.	
\end{enumerate}
\end{theorem}

\begin{proof}
(1)$\Rightarrow$(2) It follows easily from the fact that $E$ has the positive Schur property if and only if every disjointly weakly compact subset of $E$ is \textit{L}--weakly compact \cite[Corollary 2.7]{disjointly weakly compact set}.	
\par (2)$\Rightarrow$(3) Obvious.
\par (3)$\Rightarrow$(1) We assume by way of contradiction that $E$ does not have the positive Schur property. Then, there exists a disjoint  sequence $(x_n)$ of $E^+$ such that $x_{n}\xrightarrow{w}0$ and $\|x_n\|=1$ for each $n\in \mathbb{N}$, and so  we can find $ f_n\in (E^\prime)^+$ such that $\|f_n\|$=1 and $f_n(x_n)=1$ for each $n$. Hence, there exists a disjoint sequence $(g_n)\subset (E^\prime)^+$ such that
$$0\leq g_n\leq f_{n}\,\, \,  \,  \textrm{and }\,\,\, g_n(x_n)=f_n(x_n)=1.$$See, e.g., \cite[p.77]{Positive}.
 Consider the  operator $T:E\rightarrow \ell^\infty$ defined by $T(x)=(g_n(x))$ for $x\in E$. We claim that  $T$ is  \textit{LW}-compact. To this end, let $A$ be an $L$--weakly compact subset of $E$. From Proposition 3.6.2 of \cite{Meyer} it follows that $(g_n)$ converges uniformly to zero on $A$, that is, $\sup_{x\in A}|g_{n}(x)|\rightarrow0$. This implies that $T(A)\subset c_0$ and $T(A)$ is relatively compact in $c_0$, and hence $T(A)$ is relatively compact in $\ell^\infty$. So, $T$ is \textit{LW}-compact. However, $T$ is not an almost Dunford-Pettis operator since $\|T(x_n)\|\geq|g_{n}(x_n)|=1$, and hence $T$ is not a \textit{DW}-compact operator. Therefore, $E$ has the positive Schur property.
\end{proof}
 We are now in a position to consider the domination problem for positive \textit{DW}-compact operators. The following example shows that a positive operator dominated by a \textit{DW}-compact operator is not necessarily \textit{DW}-compact.
\begin{example} Let $S,T:L^{1}[0,1]\to \ell_{\infty}$ be defined by
  $$Sf=\left(\int_{0}^{1}f(t)r^{+}_{1}(t)\textrm{d}t, \int_{0}^{1}f(t)r^{+}_{2}(t)\textrm{d}t, \cdot\cdot\cdot, \int_{0}^{1}f(t)r^{+}_{n}(t)\textrm{d}t, \cdot\cdot\cdot\right)$$and
  $$Tf=\left(\int_{0}^{1}f(t)\textrm{d}t, \int_{0}^{1}f(t)(t)\textrm{d}t, \cdot\cdot\cdot, \int_{0}^{1}f(t)\textrm{d}t, \cdot\cdot\cdot\right)$$where  $(r_n)$ is the sequence of Rademacher functions on $[0,1]$. We can see that $0\leq S\leq T$ and $T$ is compact. However, $S$ is not Dunford-Pettis (hence not \textit{DW}-compact) since $r_{n}\xrightarrow{w}0$ in $L^{1}[0,1]$ and $\|Sr_n\|\geq \frac{1}{2}$.
\end{example}

\begin{theorem}\label{domination new} For two Banach lattices $E$ and $F$  the following statements are equivalent:
	\begin{enumerate}
		\item Every positive operator from $E$ into $F$ dominated by a \textit{DW}-compact operator is \textit{DW}-compact.
		\item One of the following two conditions holds:
		\begin{enumerate}
	\item [(a)] Each disjointly weakly compact subset of $E$ is $|\sigma|(E,E^{\prime})$-totally bounded, i.e.,  $E^\prime$ is discrete.
	\item [(b)] $F$ has order continuous norm.	
	\end{enumerate}
	\end{enumerate}
\end{theorem}

\begin{proof}
$(1)\Rightarrow(2)$ Assume by way of contradiction that neither is $E^\prime$ discrete nor is the norm of $F$ order continuous. Then, Aqzzouz et al. \cite{domination of AM} constructed two positive operators $S,T$ with $0\leq S\leq T$ such that  $T$ is a compact operator,  and  $S$ is not \textit{AM}-compact and hence not \textit{DW}-compact.

\par $(2)\Rightarrow(1)$ Let $S,T:E\to F$ be two positive operators such that $0\leq S\leq T$ and $T$ is \textit{DW}-compact. By Theorem \ref{AM-compact and almost DP}, $T$ is almost Dunford-Pettis and \textit{AM}-compact. Therefore, $S$ is likewise an almost Dunford-Pettis operator.
\par We assume first that the condition $(2)(a)$ holds, that is, every disjointly weakly compact subset of $E$ is $|\sigma|(E,E^{\prime})$-totally bounded. From Theorem \ref{almost DP is DW} we know that $S$ is a \textit{DW}-compact operator.
\par Now let $F$ have order continuous norm. From Theorem 4.7 of \cite{Dodds and Fremlin PL-compact} it follows that $S$ is also an \textit{AM}-compact operator since $T$ is \textit{AM}-comact (cf. e.g., \cite[Proposition 3.7.2]{Meyer}). Hence, by Theorem \ref{AM-compact and almost DP} $S$ is \textit{DW}-compact.
\end{proof}

\begin{corollary}\label{composition is DWC-compact}
  Let $E,\,F$ be two Banach lattices and let $X$ be a Banach space. Consider the scheme of operators
$E\stackrel{S_{1}}{\longrightarrow}F\stackrel{S_{2}}{\longrightarrow}X$. If the positive operator $S_{1}$ is dominated by a \textit{DW}-compact positive operator and the operator $S_{2}$ is order weakly compact, then $S_{2}{S_{1}}$ is \textit{DW}-compact.
\end{corollary}
\begin{proof} Since $S_{2}$ is order weakly compact, it follows that $S_2$ admits a factorization through  a Banach lattice $G$ with  order continuous norm

$$\xymatrix
 {E\ar[r]^{S_{1}}&F\ar[rr]^{S_{2}}\ar[dr]_{Q}&&X \\
&&G\ar[ur]_{R}&
}$$
such that $Q:F\to{G}$ is a lattice homomorphism (see, e.g., \cite[Theorem 5.58]{Positive}\,). Since the positive operator $QS_{1}:E\to{G}$ is still dominated by a \textit{DW}-compact operator and $G$ has an order continuous norm,  it follows from Theorem \ref{domination new} that $QS_{1}$ is \textit{DW}-compact. Hence $S_{2}S_{1}=RQS_{1}$ is \textit{DW}-compact.
\end{proof}

\begin{corollary}
 Consider the scheme of positive operators
$E\stackrel{S_{1}}{\longrightarrow}F\stackrel{S_{2}}{\longrightarrow}G$. If each $S_{i}$ is dominated by a \textit{DW}-compact positive operator ($i=1, 2$), then $S_{2}{S_{1}}$ is \textit{DW}-compact.
 \par In particular, if a positive operator $S$ on a Banach lattice is dominated by a \textit{DW}-compact operator, then $S^2$ is a \textit{DW}-compact operator.
\end{corollary}
\begin{proof}
 This follows from Corollary \ref{composition is DWC-compact} and the observation that $S_2$ is order weakly compact since $S_2$ is dominated by a positive \textit{DW}-compact (hence order weakly compact) operator.
\end{proof}
\vskip 5mm

\section{The (d)-DP and (d)-DP$^*$ properties of Banach lattices}

Let us recall that a bounded subset $A$ of a Banach space $X$ is called a \textit{Dunford-Pettis} (resp. \textit{limited} ) \textit{set}  if every weakly (resp. weak$^*$-) null sequence $(f_n)\subset X^\prime$ converges uniformly to zero on $A$, that is, $\sup_{x\in A}|f_{n}(x)|\to0\,\,\,(n\to \infty)$ (see \cite{Andrews, Bourgain}). We say $X$ has the \textit{Dunford-Pettis }(resp. \textit{DP}$^{\,*}$) \textit{property} whenever every relatively weakly compact set in $X$ is a Dunford-Pettis (resp. limited) set (see \cite{BFV}). Clearly, the DP$^*$ property implies the Dunford-Pettis property. It is well known that $X$ has the Dunford-Pettis property if and only if each weakly compact operator from $X$ into an arbitrary Banach space is a Dunford-Pettis operator, or equivalently, if and only if each weakly compact operator $T:X\to c_0$  is  Dunford-Pettis. Similarly, $X$ has the DP$^*$ property if and only if  each bounded linear operator $T:X\to c_0$  is  a Dunford-Pettis operator (see \cite{CGL}). Dually, a bounded subset $B$ of $X^\prime$ is called an $L$--set if every weakly null sequence $(x_n)$ of $X$ converges uniformly to zero on $B$ \cite{Em}.
\par In the past decade the disjoint and weaker versions of  the Dunford-Pettis set and limited set in Banach lattices have been introduced by Bouras \cite{Almost DP set} and Chen et al. \cite{Chen}, respectively. A bounded subset $A$ of a Banach lattice $E$ is called an \textit{almost Dunford-Pettis} (resp. \textit{almost limited} ) \textit{set}  if every disjoint weakly (resp. weak$^*$-) null sequence $(f_n)\subset E^\prime$ converges uniformly to zero on $A$. We say $E$ has the \textit{wDP }(resp. \textit{wDP}$^{\,*}$) \textit{property} whenever every relatively weakly compact set in $E$ is an almot Dunford-Pettis (resp. almost limited) set. Recently, Ardakani and Chen \cite{positively limited set} defined a class of sets called positively limited sets. A subset $A$ of $E$ is said to be a \textit{positively limited set} whenever  every weak$^*$-null sequence $(f_n)$ in $(E^\prime)^+$ converges uniformly to zero on $A$. Also, $E$ has the \textit{positive DP$^{\,*}$ property} if every relatively weakly compact subset of $E$ is positively limited.

 \par The wDP (resp. wDP$^*$, positive DP$^{*}$) property can also be characterized by using disjointly weakly compact sets instead of weakly compact sets. Indeed, a Banach lattice $E$ has the  wDP (resp. wDP$^*$, positive DP$^{*}$) property if and only if every disjointly weakly compact set in $E$ is almost Dunford-Pettis (resp. almost limited, positively limited). See \cite[Theorem 2.8 \& 2.9]{disjointly weakly compact set} and \cite[Theorem 3.10]{positively limited set}. Now we turn our attention to those Banach lattices in which each disjointly weakly compact set is a limited set.

\begin{definition}
A Banach lattice $E$ is said to have the (d)-DP$^*$ property if every disjointly weakly compact subset of $E$ is limited.
\end{definition}
$\ell^1$ is an example of a Banach lattice with the (d)-DP$^*$ property. It is obvious that the (d)-DP$^*$ property implies the DP$^*$ property. The converse does not necessarily hold. For instance, $\ell^\infty$ has the DP$^*$ property, but lacks the (d)-DP$^*$ property. Recall that Banach spaces in which the classes of relatively compact sets and limited sets coincide, are called \textit{Gelfand-Phillips spaces}. All separable Banach spaces and all WCG-spaces are Gelfand-Phillips spaces. A result due to A. V. Bukhvalov asserts that a $\sigma$-Dedekind complete Banach lattice is a Gelfand-Phillips space if and only if $E$ has order continuous norm (see, e.g.,  \cite[Theorem 4.5]{WOrder}). For our convenience,  by saying a sequence $(f_n)\subset X^\prime$ is an $L$--sequence we mean that the set $\{x_{n}:n\in \mathbb{N}\}$ is an $L$--set in $X^\prime$.
\begin{proposition}\label{disjointly DP$^*$} For a Banach lattice $E$ the following statements are equivalent.
\begin{enumerate}
  \item $E$ has the (d)-DP$^{\,*}$ property.
  \item For every weak$^*$-null sequence $(f_n)$ in $E^\prime$, we have $(f_n)$ is an $L$--sequence and $|f_n|\xrightarrow{w^*}0$.
  \item Every bounded linear operator from $E$ into a Gelfand-Phillips space is \textit{DW}-compact.
  \item Every bounded linear operator $T:E\rightarrow c_0$ is \textit{DW}-compact.
\end{enumerate}
\end{proposition}

\begin{proof}
(1) $\Rightarrow$ (2) We assume that $E$ has the (d)-DP$^{*}$ property. Let $(f_n)$ be a weak$^*$-null sequence in $E^\prime$. Since every order interval is disjointly weakly compact, hence by our assumption, is limited, it follows that $|f_n|\xrightarrow{w^*}0$. To prove that $(f_n)$ is an $L$--sequence, it suffices to show that $\lim_{n}f_{n}(x_n)=0$ for every weakly null sequence $(x_n)$ in $E$. Since every relatively weakly compact set is disjointly weakly compact,  $(x_n)$ is, again by our assumption, a limited  set. Therefore, $|f_{n}(x_n)|\leq \sup_{k}|f_{n}(x_k)|\xrightarrow{n\rightarrow\infty}0$, as desired.
\par (2) $\Rightarrow$ (1) Let $A$ be an arbitrary disjointly weakly compact subset of $E$. Without loss of generality we can assume that $A$ is solid. Let $(f_n)$ be a weak$^*$-null sequence in $E^\prime$. Then, by our hypothesis, $(f_n)$ is an $L$--sequence and $|f_n|\xrightarrow{w^*}0$. Thus, $(|f_n|)$ is also an $L$--sequence. Then, from Proposition 2.3.4 of \cite{Meyer} it follows that there exists a disjoint sequence $(x_n)\subset A^+$ such that
   $$\limsup_{n\rightarrow\infty}\rho_{A}(f_n)=\limsup_{n\rightarrow\infty}\rho_{A}(|f_n|)=\limsup_{n\rightarrow\infty}\langle|f_n|, x_n\rangle$$where $\rho_{A}(f):=\sup_{x\in A}\langle|f|,|x|\rangle$ for $f\in E^\prime$. Note that $A$ is disjointly weakly compact. It follows that $x_{n}\xrightarrow{w}0$. Since $(|f_n|)$ is also an $L$--sequence, we have
      $$ \limsup_{n\rightarrow\infty}\rho_{A}(f_n)=\limsup_{n\rightarrow\infty}\langle|f_n|, x_n\rangle=0.$$This implies that $A$ is a limited set.
\par $(1)\Rightarrow(3)\Rightarrow(4)$ are obvious since the images of limited sets under bounded linear operators are likewise limited.
\par $(4)\Rightarrow(1)$ Let $A$ be a disjointly weakly compact subset of $E$ and let $(f_n)\subset E^\prime$ be a sequence satisfying $f_{n}\xrightarrow{w^*}0$. Consider the bounded linear operator $T:E\to c_0$  defined by $Tx=(f_{n}(x))_{1}^{\infty}$ for $x\in E$. By our hypothesis $T(A)$ is relatively compact in $c_0$, that is, $\sup_{x\in A}|f_{n}(x)|\to0$. This implies that $(f_n)$ converges uniformly to zero on $A$. So, $A$ is a limited set, and hence $E$ has the (d)-DP$^{\,*}$ property.
\end{proof}
The following result shows that a $\sigma$-Dedekind complete Banach lattice with the (d)-DP$^{\,*}$ property has the Schur property.
\begin{theorem}
For a Banach lattice $E$ the following statements are equivalent:
\begin{enumerate}
  \item $E$ has the Schur property.
 \item  $E$ has the (d)-DP$^{\,*}$ property and order continuous norm.
  \item $E$ is $\sigma$-Dedekind complete  with the (d)-DP$^{\,*}$ property.
 \end{enumerate}
\end{theorem}

\begin{proof}
$(1)\Rightarrow(2)$ If $E$ has the Schur property, then it is obvious that $E$ has order continuous norm.  From Theorem \ref{Each linear operator is DW} it follows that every disjointly weakly compact subset of $E$  is relatively compact (and hence limited). This implies that $E$ has the (d)-DP$^*$ property.
\par $(2)\Rightarrow(1)$ Assume that  $E$ has the (d)-DP$^*$ property and order continuous norm. Let $A$ be a relatively weakly compact subset of $E$. Then the (d)-DP$^*$ property of $E$ implies that $A$ is limited. Since $E$ has order continuous norm, $E$ is a  Gelfand-Phillips space. Therefore, $A$ is relatively compact. This implies  that $E$ has the Schur property.
\par $(2)\Rightarrow(3)$ Obvious.
\par $(3)\Rightarrow(2)$  It suffices to show that $E$ has order continuous norm. Assume by way of contradiction that the norm of $E$ is not order continuous. Then $\ell^\infty$ lattice embeds in $E$. Since $\ell^\infty$ is an injective Banach lattice, $\ell^\infty$ can be identified with the range of a positive contradictive pojection on $E$. Since the unit ball $B_{\ell^\infty}=[-\mathbf{1},\mathbf{1}]$ is order bounded and disjointly weakly compact in $E$, by hypothesis, $B_{\ell^\infty}$ is limited in $E$ and hence $B_{\ell^\infty}$ is also limited in ${\ell^\infty}$. This is absurd.
\end{proof}

 Theorem 2.9 of \cite{disjointly weakly compact set} asserts that every disjointly weakly compact subset of $E$ is an almost  Dunford-Pettis set if and only if $E$ has the wDP property. We are now in a position to consider Banach lattices possessing the property that disjointly weakly compact sets are  Dundord-Pettis.
\begin{definition}
A Banach lattice $E$ is said to have the (d)-DP property whenever every disjointly weakly compact subset of $E$ is a Dundord-Pettis set.
\end{definition}
\begin{remark}
(1) Clearly,  the (d)-DP property implies the Dunford-Pettis property. For the converse, $\ell^{\infty}, L^{\infty}[a,b]$ and $C[a,b]$ do not have the (d)-DP property whereas they  have the Dunford-Pettis property. For KB-spaces, the (d)-DP property and the Dunford-Pettis property coincide. This can easily be seen from Proposition 2.6 of \cite{disjointly weakly compact set}. However,  $c_0$, $c$  have the (d)-DP property, but they are not KB-spaces.
\par (2) A reflexive Banach lattice $E$ has the (d)-DP property if and only if $E$ is finite dimensional.
\par (3) A well known result due to A. Grothendieck asserts  that if the dual $X^\prime$ of a Banach space $X$ has the Dunford-Pettis property, then $X$ itself has the Dunford-Pettis property. However, this is not necessarily true for the (d)-DP property. For instance, $\ell^\infty$ lacks the (d)-DP property while $(\ell^{\infty})^{\,\prime}$ is an $AL$-space and hence has the (d)-DP property.
\end{remark}

\begin{theorem}\label{d-DP}
For a Banach lattice $E$ the following statements are equivalent.
\begin{enumerate}
  \item $E$ has the (d)-DP property.
  \item $E$ has the Dunford-Pettis property and every order interval of $E$ is a Dunford-Pettis set
  \item Every weakly compact operator from $E$ into an arbitrary Banach space is a \textit{DW}-compact operator.
  \item Every weakly compact operator from $E$ into $c_0$ is a \textit{DW}-compact operator.
\end{enumerate}
\end{theorem}

\begin{proof}
\par $(1)\Rightarrow (2)$ It follows easily from a simple observation that every relatively weakly compact subset of $E$ is disjointly weakly compact and every order interval of $E$ is disjointly weakly compact.
\par $(2)\Rightarrow (3)$ Let $T:E\to X$ be a weakly compact operator from $E$ to a Banach space $X$. Then $T$ is  Dunford-Pettis since $E$ has the Dunford-Pettis property. Also, for each $x\in E^+$,  by our hypothesis, $[-x,x]$ is a Dunford-Pettis set. Therefore, $T[-x,x]$ is relatively compact. This implies that $T$ is \textit{AM}-compact. From Theorem \ref{AM-compact and almost DP} it follows that $T$ is \textit{DW}-compact.
\par $(3)\Rightarrow (4)$ Obvious.
\par $(4)\Rightarrow (1)$ Let $A$ be a disjointly weakly compact subset of $E$ and let $(f_n)\subset  E^\prime$ satisfying $f_{n}\xrightarrow{w}0$. Then the operator $T:E\to c_0$ defined by $$Tx=(f_{1}(x),f_{2}(x),\cdots,f_{n}(x),\cdots)$$is a weakly compact operator. By our hypothesis, $T$ is \textit{DW}-compact and hence $T(A)$ is a relatively compact subset of $c_0$, that is, $\sup_{x\in A}|f_{n}(x)|\to 0$. This implies that $A$ is a Dunford-Pettis set.
\end{proof}
Recall that a Banach space $X$ is said to have the \textit{reciprocal Dunford-Pettis property }if every Dunford-Pettis operator from $X$ into an arbitrary Banach space is weakly compact. A result due to C. P. Niculescu  implies that a Banach lattice $E$ has the reciprocal Dunford-Pettis property if and only if $E^{\,\prime}$ has order continuous norm \cite[Theorem 2.1]{Niculescu}. A \textit{DW}-compact operator is not necessarily weakly compact.
Indeed, we can easily see that the identity operator $I:\ell^1\to \ell^1$ is an example of a \textit{DW}-compact operator which is not weakly compact. Thus, we give the following definition.
\begin{definition}
A Banach lattice $E$ has the reciprocal (d)-DP property if every \textit{DW}-compact operator from $E$ into an arbitrary Banach space is weakly compact.
\end{definition}

The following result shows that the reciprocal (d)-DP property and the reciprocal Dunford-Pettis property indeed coincide for a Banach lattice.
\begin{proposition}
	For a Banach lattice $E$ the following statements are equivalent:
\begin{enumerate}
  \item $E$ has the reciprocal (d)-DP property.
  \item $E^{\,\prime}$ has order continuous norm.
  \item  Every \textit{DW}-compact operator from $E$ into an arbitrary Banach space is a compact operator.
  \item $E$ has the reciprocal Dunford-Pettis property
\end{enumerate}	
\end{proposition}

\begin{proof}
$(2)\Leftrightarrow(3)$ follows from Corollary \ref{DW-compact is compact}, and $(2)\Leftrightarrow(4)$ is essentially due to C. P. Niculescu (cf.  \cite[Theorem 3.7.10]{Meyer}).
\par$(3)\Rightarrow(1)$ is obvious.
 \par  $(1)\Rightarrow(2)$ Let  $T:E\rightarrow \ell^1$ be a positive operator. Then $T$ is clearly Dunford-Pettis. Since $\ell^1$ is a discrete Banach lattice with order continuous norm, $T$ is also \textit{AM}-compact. Therefore, $T$ is a \textit{DW}-compact operator. Then, by our hypothesis $T$ is weakly compact (and hence compact). In view of  a result due to B. K\"{u}hn,   $E^{\,\prime}$ has order continuous norm (cf. e.g., \cite[Theorm 5.29]{Positive}).
\end{proof}
For two fixed spaces we have the following result.
\begin{corollary} For a Banach lattice $E$ and a Banach space $X$ the following assertions are equivalent.
\begin{enumerate}
  \item Every \textit{DW}-compact operator $T:E\to X$ is a weakly compact operator.
  \item Either the norm on $E^\prime$ is order continuous  or $X$ is reflexive.
\end{enumerate}
\end{corollary}
\begin{proof}
It suffices to prove $(1)\Rightarrow(2)$. Assume by way of contradiction that neither is the norm on $E^\prime$  order continuous nor is $X$  reflexive. To finish the proof, we have to construct a \textit{DW}-compact operator $T:E\to X$ which is not weakly compact.
\par Since the norm on $E^\prime$ is not order continuous, $E$ contains a lattice isomorphic copy $U$ of $\ell^1$ and furthermore, there exists a positive projection $P:E\to U$ \cite[Proposition 2.3.11, Theorem 2.4.14]{Meyer}. Let $j:U\to\ell^1$ be the lattice isomorphism. Clearly, $jP:E\to\ell^1$ is a \textit{DW}-compact. Also, since $X$ is not reflexive, there exists a bounded sequence $(x_n)\subset X$ containing no weakly convergent subsequence. Let the operator $Q:\ell^{1}\to X$ be defined by $Q((\lambda_n))=\sum_{n=1}^{\infty}\lambda_{n}x_{n},\,\,\, (\lambda_n)\in \ell^1$. So, $QjP:E\to X$ is \textit{DW}-compact. Let $(e_n)$ be the usual unit basis of $\ell^1$. Since $QjP(j^{-1}e_{n})=x_n$, it follows that $QjP:E\to X$ is not weakly compact.
\end{proof}

\end{document}